\documentclass[reqno,10pt]{amsart}
\usepackage{amsmath,amssymb}

\usepackage[dvips]{graphics}
\usepackage{epsfig}

\newtheorem{thm}{Theorem}[section]
\newtheorem{lem}[thm]{Lemma}

\theoremstyle{definition}
\newtheorem{rek}[thm]{Remark}

\newcommand\ben{\begin{enumerate}}
\newcommand\een{\end{enumerate}}

\renewcommand{\d}{{\mathrm{d}}} 

\newcommand{\twocase}[5]{#1 \begin{cases} #2 & \text{#3}\\ #4
&\text{#5} \end{cases}  }

\newcommand\be{\begin{equation}}
\newcommand\ee{\end{equation}}
\newcommand\bea{\begin{eqnarray}}
\newcommand\eea{\end{eqnarray}}

\newcommand{\E}{{\mathbb E}} 

\renewcommand{\d}{{\mathrm{d}}} 

\numberwithin{equation}{section}

\begin{document}

\title{When the Cram\'{e}r-Rao Inequality provides no information}

\author{Steven J. Miller}
\address{Department of Mathematics, Brown University, 151 Thayer
 Street, Providence, RI 02912}
 \email{sjmiller@math.brown.edu}



\date{\today}

\begin{abstract} We investigate a one-parameter family of
probability densities (related to the Pareto distribution, which
describes many natural phenomena) where the Cram\'{e}r-Rao
inequality provides no information.
\end{abstract}

\thanks{The author would like to thank Alan Landman for many enlightening
conversations and the referees for helpful comments. The author was
partly supported by NSF grant DMS0600848.}

\maketitle


\section{Cram\'{e}r-Rao Inequality}

One of the most important problems in statistics is estimating a
population parameter from a finite sample. As there are often many
different estimators, it is desirable to be able to compare them and
say in what sense one estimator is better than another. One common
approach is to take the unbiased estimator with smaller variance.
For example, if $X_1, \dots, X_n$ are independent random variables
uniformly distributed on $[0,\theta]$, $Y_n = \max_i X_i$ and
$\overline{X} = (X_1+\cdots+X_n)/n$, then $\frac{n+1}{n} Y_n$ and $2
\overline{X}$ are both unbiased estimators of $\theta$ but the
former has smaller variance than the latter and therefore provides a
tighter estimate.

Two natural questions are (1) which estimator has the minimum
variance, and (2) what bounds are available on the variance of an
unbiased estimator? The first question is very hard to solve in
general. Progress towards its solution is given by the
Cram\'{e}r-Rao inequality, which provides a lower bound for the
variance of an unbiased estimator (and thus if we find an estimator
that achieves this, we can conclude that we have a minimum variance
unbiased estimator).

\noindent \textbf{Cram\'{e}r-Rao Inequality}: \emph{Let
$f(x;\theta)$ be a probability density function with continuous
parameter $\theta$. Let $X_1, \dots, X_n$ be independent random
variables with density $f(x;\theta)$, and let
$\widehat{\Theta}(X_1,\dots,X_n)$ be an unbiased estimator of
$\theta$. Assume that $f(x;\theta)$ satisfies two conditions:}
\begin{enumerate}

\item \emph{we have} \be \frac{\partial}{\partial \theta} \left[\int\cdots\int
\widehat{\Theta}(x_1,\dots,x_n) \prod_{i=1}^n f(x_i;\theta) \d x_i
\right] \ = \ \int\cdots\int \widehat{\Theta}(x_1,\dots,x_n)
\frac{\partial \prod_{i=1}^n f(x_i;\theta)}{\partial \theta} \d
x_1\cdots \d x_n; \ee

\item \emph{for each $\theta$, the variance of
$\widehat{\Theta}(X_1,\dots,X_n)$ is finite.}
\end{enumerate}

\emph{Then} \be\label{eq:crineq} {\rm var}(\widehat{\Theta}) \ \ge \
\cfrac1{n \E\left[ \left( \frac{\partial \log f(x; \theta)}{\partial
\theta}\right)^2\right]}, \ee \emph{where $\E$ denotes the expected
value with respect to the probability density function
$f(x;\theta)$.}\\

For a proof, see for example \cite{CaBe}. The expected value in
\eqref{eq:crineq} is called the \emph{information number} or the
\emph{Fisher information} of the sample.

As variances are non-negative, the Cram\'{e}r-Rao inequality
(equation \eqref{eq:crineq}) provides no useful bounds on the
variance of an unbiased estimator if the information is infinite, as
in this case we obtain the trivial bound that the variance is
greater than or equal to zero. We find a simple one-parameter family
of probability density functions (related to the Pareto
distribution) that satisfy the conditions of the Cram\'{e}r-Rao
inequality, but the expectation (i.e., the
information) is infinite. Explicitly, our main result is\\

\noindent \textbf{Theorem:} \emph{Let}\be \twocase{f(x;\theta) \ = \
}{a_\theta\ x^{-\theta}\ \log^{-3} x}{if $x \ge e$}{0}{otherwise,}
\ee \emph{where $a_\theta$ is chosen so that $f(x;\theta)$ is a
probability density function. The information is infinite when
$\theta = 1$. Equivalently, the Cram\'{e}r-Rao inequality yields the
trivial (and useless) bound that ${\rm Var}(\widehat{\Theta}) \ge 0$
for any unbiased estimator $\widehat{\Theta}$ of $\theta$ when
$\theta = 1$.}

\bigskip

In \S\ref{sec:almostpareto} we analyze the density in our theorem in
great detail, deriving needed results about $a_\theta$ and its
derivatives as well as discussing how $f(x;\theta)$ is related to
important distributions used to model many natural phenomena. We
show the information is infinite when $\theta=1$ in
\S\ref{sec:computeinfo}, which proves our theorem. We also discuss
there properties of estimators for $\theta$. While it is not clear
whether or not this distribution has an unbiased estimator, there is
(at least for $\theta$ close to 1) an asymptotically unbiased
estimator rapidly converging to $\theta$ as the sample size tends to
infinity. By examining the proof of the Cram\'{e}r-Rao inequality we
see that we may weaken the assumption of an unbiased estimator.
While typically there is a cost in such a generalization, as our
information is infinite there is no cost in our case. We may
therefore conclude that arguments such as those used to prove the
Cram\'{e}r-Rao inequality cannot provide any information for
estimators of $\theta$ from this distribution.


\section{An Almost Pareto Density}\label{sec:almostpareto}

Consider \be \twocase{f(x;\theta) \ = \ }{a_\theta / (x^\theta
\log^3 x)}{if $x \ge e$}{0}{otherwise,} \ee where $a_\theta$ is
chosen so that $f(x;\theta)$ is a probability density function. Thus
\be \int_{e}^\infty a_\theta \frac{\d x}{x^\theta \log^3 x} \ = \ 1.
\ee We chose to have $\log^3 x$ in the denominator to ensure that
the above integral converges, as does $\log x$ times the integrand;
however, the expected value (in the expectation in
\eqref{eq:crineq}) will not converge.

For example, $1 / x \log x$ diverges (its integral looks like
$\log\log x$) but $1/x\log^2 x$ converges (its integral looks like
$1/\log x$); see pages 62--63 of \cite{Rud} for more on close
sequences where one converges but the other does not. This
distribution is close to the Pareto distribution (or a power law).
Pareto distributions are very useful in describing many natural
phenomena; see for example \cite{DM,Ne,NM}. The inclusion of the
factor of $\log^{-3} x$ allows us to have the exponent of $x$ in the
density function equal $1$ \emph{and have the density function
defined for arbitrarily large $x$}; it is also needed in order to
apply the Dominated Convergence Theorem to justify some of the
arguments below. If we remove the logarithmic factors then we obtain
a probability distribution only if the density vanishes for large
$x$. As $\log^3 x$ is a very slowly varying function, our
distribution $f(x;\theta)$ may be of use in modeling data from an
unbounded distribution where one wants to allow a power law with
exponent $1$, but cannot as the resulting probability integral would
diverge. Such a situation occurs frequently in the Benford Law
literature; see \cite{Hi,Rai} for more details.

We study the variance bounds for unbiased estimators
$\widehat{\Theta}$ of $\theta$, and in particular we show that when
$\theta = 1$ then the Cram\'{e}r-Rao inequality yields a useless
bound.

Note that it is not uncommon for the variance of an unbiased
estimator to depend on the value of the parameter being estimated.
For example, consider again the uniform distribution on
$[0,\theta]$. Let $\overline{X}$ denote the sample mean of $n$
independent observations, and $Y_n = \max_{1 \le i \le n} X_i$ be
the largest observation. The expected value of $2\overline{X}$ and
$\frac{n+1}{n} Y_n$ are both $\theta$ (implying each is an unbiased
estimator for $\theta$); however, ${\rm Var}(2\overline{X}) =
\theta^2/3n$ and ${\rm Var}(\frac{n+1}{n}Y_n) = \theta^2/n(n+1)$
both depend on $\theta$, the parameter being estimated (see, for
example, page 324 of \cite{MM} for these calculations).

\begin{lem}\label{lem:athetastuff}
As a function of $\theta \in [1,\infty)$, $a_\theta$ is a strictly
increasing function and $a_1 = 2$. It has a one-sided derivative at
$\theta = 1$, and $\frac{\d a_\theta}{\d \theta} \in (0,\infty)$.
\end{lem}

\begin{proof} We have \be a_\theta \int_e^\infty \frac{\d x}{x^\theta \log^3
x} \ = \ 1. \ee When $\theta = 1$ we have \be\label{eq:athetaint}
a_1 \ = \ \left[ \int_e^\infty \frac{\d x}{x\log^3 x} \right]^{-1},
\ee which is clearly positive and finite. In fact, $a_1 = 2$ because
the integral is \be \int_e^\infty \frac{\d x}{x \log^3 x} \ = \
\int_e^\infty \log^{-3}x\ \frac{\d\log x}{\d x} \ = \
\frac{-1}{2\log^2 x}\Bigg|_e^\infty \ = \ \frac12; \ee though all we
need below is that $a_1$ is finite and non-zero, we have chosen to
start integrating at $e$ to make $a_1$ easy to compute.

It is clear that $a_\theta$ is strictly increasing with $\theta$, as
the integral in \eqref{eq:athetaint} is strictly decreasing with
increasing $\theta$ (because the integrand is decreasing with
increasing $\theta$).

We are left with determining the one-sided derivative of $a_\theta$
at $\theta=1$, as the derivative at any other point is handled
similarly (but with easier convergence arguments). It is technically
easier to study the derivative of $1/a_\theta$, as
\be\label{eq:datdaoot} \frac{\d}{\d \theta} \frac1{a_\theta} \ = \
-\frac1{a_\theta^2} \frac{\d a_\theta}{\d \theta} \ee and \be
\frac1{a_\theta} \ = \ \int_e^\infty \frac{\d x}{x^\theta \log^3 x}.
\ee The reason we consider the derivative of $1/a_\theta$ is that
this avoids having to take the derivative of the reciprocals of
integrals. As $a_1$ is finite and non-zero, it is easy to pass to
$\frac{\d a_\theta}{\d \theta}|_{\theta=1}$. Thus we have \bea
\frac{\d}{\d \theta} \frac1{a_\theta}\Big|_{\theta = 1} & \ = \ &
\lim_{h\to 0^+} \frac1{h}\left[ \int_e^\infty \frac{\d
x}{x^{1+h}\log^3 x} - \int_e^\infty \frac{\d x}{x\log^3 x}\right]
\nonumber\\ & = & \lim_{h\to 0^+} \int_e^\infty \frac{1 - x^h}{h}
\frac{1}{x^h} \frac{\d x}{x\log^3 x}. \eea We want to interchange
the integration with respect to $x$ and the limit with respect to
$h$ above. This interchange is permissible by the Dominated
Convergence Theorem (see Appendix \ref{sec:applydct} for details of
the justification). Note \be\label{eq:lhospital} \lim_{h\to 0^+}
\frac{1-x^h}{h} \frac{1}{x^h} \ = \ -\log x; \ee one way to see this
is to use the limit of a product is the product of the limits, and
then use L'Hospital's rule, writing $x^h$ as $e^{h\log x}$.
Therefore \be\label{eq:dooathetaint} \frac{\d}{\d \theta}
\frac1{a_\theta}\Big|_{\theta=1} \ = \ - \int_e^\infty \frac{\d x}{x
\log^2 x}; \ee as this is finite and non-zero, this completes the
proof and shows $\frac{\d a_\theta}{\d \theta}|_{\theta=1} \in
(0,\infty)$.
\end{proof}

\begin{rek} We see now why we chose $f(x;\theta) = a_\theta
/x^\theta\log^3 x$ instead of $f(x;\theta) = a_\theta
/x^\theta\log^2 x$. If we only had two factors of $\log x$ in the
denominator, then the one-sided derivative of $a_\theta$ at
$\theta=1$ would be infinite.
\end{rek}

\begin{rek}\label{rek:derivatheta1} Though the actual value of
$\frac{\d a_\theta}{\d \theta}|_{\theta=1}$ does not matter, we can
compute it quite easily. By \eqref{eq:dooathetaint} we have \bea
\frac{\d}{\d \theta} \frac1{a_\theta}\Big|_{\theta=1}
& \ = \ & -\int_e^\infty \frac{\d x}{x\log^2 x} \nonumber\\
& = & -\int_e^\infty \log^{-2}x\ \frac{\d\log x}{\d x} \nonumber\\ &
= & \frac{1}{\log x}\Big|_e^\infty \ = \ -1. \eea Thus by
\eqref{eq:datdaoot}, and the fact that $a_1 = 2$ (Lemma
\ref{lem:athetastuff}), we have \be \frac{\d a_\theta}{\d
\theta}\Big|_{\theta=1} \ = \ -a_1^2 \cdot \frac{\d}{\d \theta}
\frac1{a_\theta}\Big|_{\theta=1} \ = \ 4. \ee
\end{rek}


\section{Computing the Information}\label{sec:computeinfo}

We now compute the expected value, $\E\left[ \left( \frac{\partial
\log f(x; \theta)}{\partial \theta}\right)^2\right]$; showing it is
infinite when $\theta=1$ completes the proof of our main result.
Note \bea \log f(x;\theta) & \ = \ & \log
a_\theta -\theta \log x + \log \log^{-3} x \nonumber\\
\frac{\partial \log f(x;\theta)}{\partial \theta} & = &
\frac{1}{a_\theta} \frac{\d a_\theta}{\d \theta} - \log x. \eea

By Lemma \ref{lem:athetastuff} we know that $\frac{\d a_\theta}{\d
\theta}$ is finite for each $\theta \ge 1$. Thus \bea \E\left[
\left( \frac{\partial \log f(x; \theta)}{\partial
\theta}\right)^2\right] & \ = \ & \E\left[ \left( \frac{1}{a_\theta}
\frac{\d a_\theta}{\d \theta} - \log x \right)^2\right] \nonumber\\
& = & \int_e^\infty \left(\frac{1}{a_\theta} \frac{\d a_\theta}{\d
\theta} - \log x \right)^2 \cdot a_\theta \frac{\d x}{x^\theta
\log^3 x}. \eea If $\theta > 1$ then the expectation is finite and
non-zero. We are left with the interesting case when $\theta = 1$.
As $\frac{\d a_\theta}{\d \theta}|_{\theta=1}$ is finite and
non-zero, for $x$ sufficiently large (say $x \ge x_1$ for some
$x_1$, though by Remark \ref{rek:derivatheta1} we see that we may
take any $x_1 \ge e^4$) we have \be \left|\frac{1}{a_1} \frac{\d
a_\theta}{\d \theta}\Big|_{\theta=1} \right| \ \le \ \frac{\log x}2.
\ee As $a_1 = 2$, we have \bea \E\left[ \left( \frac{\partial \log
f(x; \theta)}{\partial \theta}\right)^2\right]\Bigg|_{\theta=1} & \
\ge \ & \int_{x_1}^\infty
\left(\frac{\log x}2\right)^2 a_1 \frac{\d x}{x\log^3 x} \nonumber\\
& = & \int_{x_1}^\infty \frac{\d x}{2x\log x} \nonumber\\  & \ = \ &
\frac12\int_{x_1}^\infty \log^{-1} x\ \frac{\d\log x}{\d x}
\nonumber\\ & = & \frac12 \log\log x\Big|_{x_1}^\infty \nonumber\\
& = & \infty. \eea Thus the expectation is infinite. Let
$\widehat{\Theta}$ be \emph{any} unbiased estimator of $\theta$. If
$\theta = 1$ then the Cram\'{e}r-Rao inequality gives \be {\rm
var}(\widehat{\Theta}) \ \ge \ 0, \ee which provides no information
as variances are always non-negative. This completes the proof of
our theorem. \hfill $\Box$

\begin{figure}
\begin{center}
\scalebox{.5}{\includegraphics{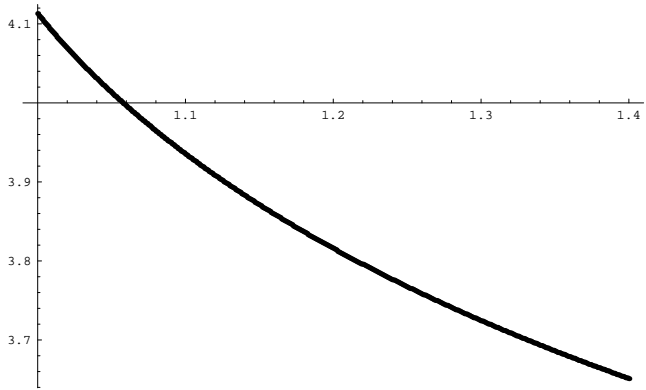}}
\caption{\label{fig:medianvtheta}Plot of the median
$\tilde{\mu}_\theta$ of $f(x;\theta)$ as a function of $\theta$
($\tilde{\mu}_1 = e^{\sqrt{2}}$).}
\end{center}\end{figure}

We now discuss estimators for $\theta$ for our distribution
$f(x;\theta)$. If $X_1, \dots, X_n$ are $n$ independent random
variables with common distribution $f(x;\theta)$, then as
$n\to\infty$ the sample median converges to the population median
$\tilde{\mu}_\theta$ (if $n=2m+1$ then the sample median converges
to being normally distributed with median $\widetilde{\mu}_\theta$
and variance $1/8mf(\tilde{\mu}_\theta;\theta)^2$; see for example
Theorem 8.17 of \cite{MM}). For $\theta$ close to 1 we see in Figure
\ref{fig:medianvtheta} that the median $\tilde{\mu}_\theta$ of
$f(x;\theta)$ is strictly decreasing with increasing $\theta$,
which implies that there is an inverse function $g$ such that
$g(\tilde{\mu}_\theta) = \theta$. We obtain an estimator to $\theta$
by applying $g$ to the sample median. This estimator is a consistent
estimator (as the sample size tends to infinity it will tend to
$\theta$) and should be asymptotically unbiased.

The proof of the Cram\'{e}r-Rao inequality starts with \be 0 \ = \
\E\left[\int \cdots \int \left(\widehat{\Theta}(x_1,\dots,x_n) -
\theta\right) h(x_1;\theta) \cdots h(x_n;\theta)dx_1\cdots
dx_n\right], \ee where $\widehat{\Theta}(x_1,\dots,x_n)$ is  an
unbiased estimator of $\theta$ depending only on the sample values
$x_1, \dots, x_n$. In our case (when each $h(x;\theta) =
f(x;\theta)$) we may not have an unbiased estimator. If we denote
this expectation by $\mathcal{F}(\theta)$, for our investigations
all that we require is that $d\mathcal{F}(\theta)/d\theta$ is finite
(which is easy to show). Going through the proof of the
Cram\'{e}r-Rao inequality shows that the effect of this is to
replace the factor of $1$ in \eqref{eq:crineq} with $(1 +
d\mathcal{F}(\theta)/d\theta)^2$; thus the generalization of the
Cram\'{e}r-Rao inequality for our estimator is \be {\rm
var}(\widehat{\Theta}) \ \ge \
\left(1+\frac{d\mathcal{F}(\theta)}{d\theta}\right)^2 \ \ \Big/\ \ n
\E\left[ \left( \frac{\partial \log f(x; \theta)}{\partial
\theta}\right)^2\right]. \ee As our variance is infinite for
$\theta=1$ we see that, no matter what `nice' estimator we use, we
will not obtain any useful information from such arguments.




\appendix

\section{Applying the Dominated Convergence
Theorem}\label{sec:applydct}

We justify applying the Dominated Convergence Theorem in the proof
of Lemma \ref{lem:athetastuff}. See, for example, \cite{SS} for the
conditions and a proof of the Dominated Convergence Theorem.

\begin{lem} For each fixed $h > 0$ and any $x \ge e$, we have
\be\label{eq:dctneededest} \left|\frac{1-x^h}{h}
\frac{1}{x^h}\right|\ \ \le\ \ e \log x,\ee and $\frac{e\log
x}{x\log^3 x}$ is positive and integrable, and dominates each
$\frac{1-x^h}{h} \frac{1}{x^h} \frac{1}{x\log^3 x}$.
\end{lem}

\begin{proof} We first prove \eqref{eq:dctneededest}. As $x \ge e$ and $h > 0$,
note $x^h \ge 1$. Consider the case of $1/h \le \log x$. Since
$|1-x^h| < 1+x^h \le 2 x^h$, we have \be \frac{|1-x^h|}{h x^h}\ \ <
\ \ \frac{2x^h}{hx^h}\ \ \le\ \ \frac2h\ \ \le\ \ 2 \log x. \ee We
are left with the case of $1/h > \log x$, or $h\log x < 1$. We have
\bea |1 - x^h| &\ \ = \ \ & |1 - e^{h\log x}| \nonumber\\ & = &
\left|1 - \sum_{n=0}^\infty \frac{(h\log x)^n}{n!}\right|
\nonumber\\ & = & h\log x \sum_{n=1}^\infty \frac{(h\log
x)^{n-1}}{n!} \nonumber\\ & < & h\log x \sum_{n=1}^\infty
\frac{(h\log x)^{n-1}}{(n-1)!}\ \ = \ \ h \log x \cdot e^{h\log x}.
\eea This, combined with $h\log x < 1$ and $x^h \ge 1$ yields \bea
\frac{|1-x^h|}{h x^h}\ \ <\ \ \frac{eh\log x}{h}\ \ = \ \ e \log x.
\eea It is clear that $\frac{\log x}{x\log^3 x}$ is positive and
integrable, and by L'Hospital's rule (see \eqref{eq:lhospital}) we
have that \be \lim_{h\to 0^+} \frac{1-x^h}{h} \frac{1}{x^h}
\frac{1}{x\log^3 x}\ \ = \ \ -\frac1{x\log^2 x}. \ee Thus the
Dominated Convergence Theorem implies that \be \lim_{h\to 0^+}
\int_e^\infty \frac{1-x^h}{h} \frac{1}{x^h} \frac{\d x}{x\log^3 x} \
\ = \ \ -\int_e^\infty \frac{\d x}{x \log^2 x} \ \ = \ \ -1 \ee (the
last equality is derived in Remark \ref{rek:derivatheta1}).
\end{proof}



\bigskip

\end{document}